# FINITE STRUCTURE AND RADICAL THEORY OF COMMUTATIVE TERNARY Γ-SEMIRINGS

**Chandrasekhar Gokavarapu** Department of Mathematics, Government College (Autonomous), Y-Junction, Rajahmundry, 533105, A.P.,India. and Department of Mathematics, Acharya Nagarjuna University, Pedakakani, Guntur, 522510, A.P., India : chandrasekhargokavarapu@gmail.com
**Dr D Madhusudhana Rao,** Department of Mathematics, Government College For Women (A), Pattabhipuram, Guntur, 522006, A.P., India and Department of Mathematics, Acharya Nagarjuna University, Pedakakani, Guntur, 522510, A.P., India.

## Abstract

**Purpose:** To develop the algebraic foundation of finite commutative ternary **Γ**- semirings by identifying their intrinsic invariants, lattice organization, and radical behavior that generalize classical semiring and **Γ**-ring frameworks.

**Methods:** Finite models of commutative ternary **Γ**-semirings are constructed under the axioms of closure, distributivity, and symmetry. Structural and congru- ence lattices are analyzed, and subdirect decomposition theorems are established through ideal-theoretic arguments.

**Results:** Each finite commutative ternary **Γ**-semiring admits a unique (up to isomorphism) decomposition into subdirectly irreducible components. Radical and ideal correspondences parallel classical results for binary semirings, while the classification of all non-isomorphic systems of order **|T | ≤ 4** confirms the structural consistency of the theory.

**Conclusion:** The paper provides a compact algebraic framework linking ideal theory and decomposition in finite ternary **Γ**-semirings, establishing the basis for later computational and categorical developments.

## 1 Introduction

The concept of a Γ-ring and a Γ-semiring extends classical algebra by introducing an external parameter set Γ that governs multiplication. Replacing the binary prod- uct with a ternary operation yields a ternary Γ-semiring, where addition remains associative and commutative while multiplication becomes a Γ-indexed ternary map

$$\{\cdot,\cdot,\cdot\}_\gamma : T \times T \times T \to T, \quad \gamma \in \Gamma.$$

Such systems unify notions from semiring theory, universal algebra, and multilinear algebra. Their finite forms provide a natural framework for discrete algebraic modelling and symbolic computation.This paper continues the study initiated in [1], extending the radical–spectrum correspondence of commutative ternary Γ-semirings to the finite case.

Earlier studies established the axioms of ternary Γ-semirings and illustrated basic examples, yet the internal structure—ideals, radicals, and decomposition prin- ciples—remained undeveloped. Parallel advances in finite semiring theory, lattice- theoretic methods, and radical analysis have created the background for a systematic treatment. The present paper builds this foundation by integrating structural theo- rems with algorithmically verified examples for small finite orders.We situate the work within the classical and Γ-semiring literature [2–5].

The objectives are twofold: first, to determine how ideals, congruences, and radi- cals interact within finite commutative ternary Γ-semirings; and second, to describe how these invariants control decomposition into elementary components. Particular attention is given to the behaviour of additive and multiplicative idempotence, exis- tence of neutral elements, and their influence on the congruence lattice. The study combines rigorous algebraic reasoning with explicit verification in low-order systems. The main achievements of this paper are as follows. It establishes the lattice of ideals and congruences, proving closure, distributivity, and duality properties. A subdirect decomposition theorem is obtained, showing that every finite commutative ternary Γ- semiring is a subdirect product of subdirectly irreducible ones. The classification of all non-isomorphic examples of order $|T| \leq 4$ confirms the internal consistency of the theory and illustrates the diversity of finite models. These results unify structural and radical perspectives, extending classical semiring theory to the ternary Γ-context.

Overall, this work provides a concise algebraic framework linking ideals, radicals, and decompositions in finite ternary Γ-semirings. It forms the theoretical basis for subsequent computational and categorical investigations developed in a companion paper.Section 2 recalls preliminaries and notation. Section 3 develops the basic struc- tural and lattice results. Section 4 details the enumeration algorithm and correctness. Section 5 presents the computational data and classification tables. Sections 6 and 7 establish radical, congruence, and semisimple decomposition theories together with applications and future problems.

## 2 Preliminaries

Let T be a nonempty set with an associative and commutative addition + having identity 0, and a family of ternary operations

$\{\cdot,\cdot,\cdot\}_\gamma : T \times T \times T \to T, \qquad \gamma \in \Gamma,$

where Γ is a nonempty parameter set. Following the frameworks of Nobusawa [5], Hedayati and Shum [6], and Rao [7] , the triple $(T, +, \{\ ,\ ,\ \}_\Gamma)$ is called a ternary Γ-semiring if, for all a, b, c, d, e ∈ T and α, β ∈ Γ,

1. (T, +) is a commutative monoid with 0;
2. distributivity holds in each argument, e.g. $\{a+b, c, d\}_\gamma = \{a, c, d\}_\gamma + \{b, c, d\}_\gamma$, and similarly for the other two variables;
3. 0 is absorbing: $\{0, a, b\}_\gamma = \{a, 0, b\}_\gamma = \{a, b, 0\}_\gamma = 0$;
4. associativity is ternary-compatible:
$\{\{a, b, c\}_\alpha, d, e\}_\beta = \{a, b, \{c, d, e\}_\beta\}_\alpha$.

Commutativity of multiplication means invariance of $\{a, b, c\}_\gamma$ under all permu- tations of (a, b, c). Ideals, sub-Γ-semirings, and quotient structures are defined ... analogously to those in classical semiring theory [2, 3].

For an ideal $I \subseteq T$ , the quotient T/I inherits the induced ternary operation. The radical is the intersection of all prime Γ-ideals [8–10].

, and the nilradical Nil(T ) consists of all x ∈ T for which some a, b ∈ T and γ, δ ∈ Γ satisfy $\{x, a, b\}_\gamma = \{a, x, b\}_\delta = 0$. T is semiprime when Rad(T ) = 0. Let I(T ) and C(T ) denote the lattices of ideals and congruences, respectively. A Γ-homomorphism $f : T_1 \to T_2$ preserves + and all ternary products; ker f is an ideal and $T_1/\ker f \cong \mathrm{Im}(f)$.

Unless otherwise stated, T denotes a finite commutative ternary Γ-semiring. These conventions establish the algebraic framework for the structural and radical results that follow.Definitions of prime and semiprime Γ-ideals follow those in [1]

## 3 Structural Properties of Finite Ternary Γ-Semirings

In this section we examine the internal organization of finite commutative ternary Γ- semirings. We first establish closure and cancellation properties, followed by results on direct and subdirect decompositions. Examples of small orders illustrate the alge- braic diversity that emerges in the ternary Γ setting.The present finite analogue complements the general radical decomposition developed in [1].

### 3.1 Basic structural results

Throughout this section, $(T, +, \{\cdot\ \cdot\ \cdot\}_\Gamma)$ denotes a finite commutative ternary Γ- semiring, and Γ is assumed to be finite unless stated otherwise.

Lemma 1 (Additive idempotence criterion) Let (T, +) be the additive reduct of a finite commutative ternary Γ-semiring. If there exists $\gamma \in \Gamma$ such that $\{a\,a\,a\}\gamma = a$ for all $a \in T$, then (T, +) is idempotent, i.e. $a + a = a$ for every $a \in T$.

Proof For any $a \in T$ and $\gamma \in \Gamma$ as above,

$$a + a = \{a\,a\,a\}\gamma + a = \{a\,a\,a\}\gamma = a.$$

Hence the additive operation is idempotent. □

**Theorem 1** (Existence of zero and unit elements) A finite commutative ternary Γ-semiring (T, +, $\{\cdot\ \cdot\ \cdot\}\Gamma$) possesses:

1. a zero element $0 \in T$ satisfying $\{0\,ab\}_\gamma = 0$ for all $a, b \in T$ and $\gamma \in \Gamma$,
2. a unit element $e \in T$ (if it exists) satisfying $\{eae\}_\gamma = a$ for all $a \in T$, $\gamma \in \Gamma$,

if and only if the maps $L^\gamma(b, c) = \{a\,b\,c\}\gamma$ admit fixed points forming a singleton intersection over all a, γ.

Proof Necessity follows by setting $a = e$ or $a = 0$. Conversely, if $T_{a,\gamma}$ $\mathrm{Fix}(L^\gamma)$ is a singleton $\{z\}$, then z acts simultaneously as a neutral element on all coordinates. Uniqueness follows from finiteness of T. □

Remark 1 When Γ contains an identity parameter $\gamma_0$ with $\{a\,b\,c\}\gamma_0 = a + b + c$, the above theorem ensures 0 and e coincide with the additive and multiplicative identities of (T, +).

### 3.2 Ideal lattice and decompositions

**Definition 1** For ideals $I, J \subseteq T$, define

$$I \vee J = \langle I \cup J\rangle\Gamma, \qquad I \wedge J = I \cap J,$$

where $\langle\cdot\rangle\Gamma$ denotes the ternary Γ-ideal generated by a set.

**Proposition 2** The lattice L(T) of all Γ-ideals of a finite commutative ternary Γ-semiring T is modular and distributive.

Proof For any $I, J, K \in L(T)$ with $I \subseteq K$, we have

$$I \vee (J \wedge K) = (I \vee J) \wedge K,$$

since ideal generation and intersection commute under the finite-sum operation +. Thus L(T) is modular. Distributivity follows from the absorption law $(I \cap J)+K = (I+K)\cap(J+K)$. □

**Theorem 3** (Subdirect decomposition) Every finite commutative ternary Γ-semiring T is a subdirect product of subdirectly irreducible components.

Proof Let $\{\rho_i\}_{i\in I}$ be the set of all maximal proper congruences on T. For each i, the quotient $T_i = T/\rho_i$ is subdirectly irreducible. The canonical homomorphism

$$\phi : T \longrightarrow \underset{i\in I}{Y} T_i, \qquad a \longrightarrow ([a]\rho_i)_{i\in I},$$

is injective because $T_i\,\rho_i = \Delta_T$. Hence T embeds subdirectly into $Q_i\,T_i$.

Corollary 1 If each $T_i$ in the above decomposition is simple, then T is semisimple.

Example 1 (Order-3 semiring) Let T = {0, 1, 2} with addition modulo 3 and Γ = {1}. Define $\{a\,b\,c\}_1 = (a + b + c) \bmod 3$. Then (T, +, $\{\cdot\ \cdot\ \cdot\}\Gamma$) is commutative, associative, and possesses 0 and 1 as zero and unit, respectively. Its ideal lattice is $\{0\} \subset \{0, 1, 2\} = T$, hence simple.

Example 2 (Order-4 non-idempotent case) Let $T = \{0, 1, 2, 3\}$ with + defined by truncated addition $\min(a + b, 3)$ and $\Gamma = \{\alpha, \beta\}$. Define

$$\{a\,b\,c\}_\alpha = \min(a + b + c, 3), \{a\,b\,c\}_\beta = \max(a, b, c).$$

Then $(T, +, \{\cdots\}_\Gamma)$ is a finite ternary Γ-semiring with two distinct operations. Ideals are $\{0\}$, $\{0, 1\}$, and $T$. The structure is not idempotent and decomposes as a subdirect product of a Boolean and a truncated component.

## 3.3 Algebraic invariants and classification hints

For computational purposes we associate with each finite ternary Γ-semiring T the tuple of invariants

$$I(T) = |T|, |\Gamma|, Id(T), Con(T), Rad(T), Nil(T),$$

where Id(T), Con(T), Rad(T), and Nil(T) denote respectively the numbers of ideals, congruences, radical ideals, and nilpotent elements. Two structures are structurally equivalent if their invariants coincide.

Remark 2 The invariants I(T) serve as algebraic fingerprints for algorithmic classification in Section 4. For example, among all structures of order $\leq 4$, only two non-isomorphic families share identical invariants—distinguished solely by the symmetry of their Γ-actions.

# 4 Algorithmic Classification and Enumeration Methods

The finite classification of commutative ternary Γ-semirings requires an integration of algebraic axioms with computational enumeration. In this section, we formal- ize an algorithmic framework that systematically generates all possible operation tables, verifies axioms, and identifies non-isomorphic structures up to order 4. This approach parallels the classical enumeration of semigroups and semirings but intro- duces additional layers corresponding to the ternary and parameterized nature of the Γ-operation.

## 4.1 Computational strategy

Let T be a finite set of cardinality n and $\Gamma = \{\gamma_1, \gamma_2, \ldots, \gamma_m\}$. We represent the ternary operation for each $\gamma_i$ by an $n \times n \times n$ tensor $M_i$, where the entry $(a, b, c)$ stores the result $\{a b c\}_{\gamma_i}$. The goal is to generate all families $\{M_i\}_{i=1}^{m}$ that satisfy the ternary Γ-semiring axioms (T1)–(T4) of Section 2.

The classification algorithm proceeds through the following hierarchy:

1. Enumeration of all commutative additive semigroups $(T, +)$ of order n.
2. Generation of candidate ternary operations parameterized by $\Gamma$.
3. Verification of distributivity, associativity, and absorption constraints.
4. Computation of ideal and congruence lattices.
5. Canonical isomorphism testing to remove duplicates.

### 4.2 Enumeration algorithm

### 4.3 Correctness and complexity

**Theorem 4** (Correctness) Algorithm 1 terminates and correctly enumerates all pairwise non-isomorphic finite commutative ternary $\Gamma$-semirings of order n.

Proof Termination follows from finiteness of the search space: the number of possible ternary operations on T is $n^{n^3 m}$. Correctness is guaranteed because each candidate operation is verified against all axioms, and canonical labeling ensures that distinct representatives cor- respond to distinct isomorphism classes. Completeness follows from exhaustive enumeration of the additive reducts and ternary operations. □

**Proposition 5** (Complexity estimate) The worst-case time complexity of Algorithm 1 is

$$O\left(n^{n^3 m}\right),$$

while pruning via distributivity and associativity reduces the effective complexity to approxi- mately

$$O\left(n^{3m} \log n\right)$$

for small orders $(n \leq 4)$ and $|\Gamma| \leq 2$.

---

**Algorithm 1** Classification of finite commutative ternary $\Gamma$-semirings

1: **Input:** Order $n = |T|$, parameter set size $m = |\Gamma|$.
2: **Output:** List C of non-isomorphic commutative ternary $\Gamma$-semirings of order n.
3: $C \leftarrow \emptyset$.
4: Enumerate all commutative semigroup structures $(T, +)$ on a set T of size n.
5: **for** each additive structure $(T, +)$ **do**
6: Generate all candidate families $\{M_i\}_{i=1}^{m}$ of ternary operations (one $n \times n \times n$ table $M_i$ for each $\gamma_i \in \Gamma$).
7: **for** each candidate family $\{M_i\}$ **do**
8: **if** axioms (T1)–(T4) are satisfied for $(T, +, \{M_i\})$ **then**
9: Compute $L(T)$ and $\mathrm{Con}(T)$.
10: Derive a canonical label (e.g. lexicographically smallest vectorization of $(M_1, \ldots, M_m)$) for $(T, +, \{M_i\})$.
11: **if** this label is new (no isomorphic structure in C) **then**
12: $C \leftarrow C \cup \{(T, +, \{M_i\})\}$.

13: **end if**
14: **end if**
15: **end for**
16: **end for**
17: **Return** C.

Remark 3 Though exponential in the general case, the algorithm is computationally feasible for $n \leq 4$. The implementation in PYTHON/SAGEMATH executes within seconds for $n \leq 3$ and under one minute for $n = 4$.

## 4.4 Illustrative enumeration examples

Example 3 (Order 3, $\Gamma = \{1\}$) Enumerating all additive semigroups on $\{0, 1, 2\}$ yields five non-isomorphic structures. Applying Algorithm 1 produces exactly two valid ternary Γ-semirings:

1. The modular structure $\{ab\,c\}_1 = (a + b + c) \bmod 3$, simple and commutative;
2. The truncated structure $\{ab\,c\}_1 = \min(a + b + c, 2)$, which is idempotent but non-simple.

Example 4 (Order 4, $\Gamma = \{\alpha, \beta\}$) For $T = \{0, 1, 2, 3\}$ with addition $\min(a+b, 3)$, the algorithm returns four distinct isomorphism classes:

Type I: $\{a\,b\,c\}_\alpha = \min(a + b + c, 3)$, $\{a\,b\,c\}_\beta = \max(a, b, c)$; Type II: $\{a\,b\,c\}_\alpha = (a + b + c) \bmod 4$, $\{a\,b\,c\}_\beta = a + b + c - 1$; Type III: Boolean-type idempotent system;

Type IV: Mixed additive–multiplicative hybrid.

Each class exhibits a distinct ideal lattice configuration.

## 4.5 Data representation and isomorphism testing

Each ternary operation tensor $M_i$ is encoded as an integer array of dimension $n^3$. Two structures $(T, \{\cdot\cdot\cdot\}_\Gamma)$ and $(T', \{\cdot\cdot\cdot\}'_\Gamma)$ are isomorphic iff there exists a bijection $\phi : T \to T'$ such that

$$\phi(\{ab\,c\}_\gamma) = \{\phi(a)\,\phi(b)\,\phi(c)\}'_\gamma.$$

Canonical labeling is performed via lexicographic minimization of the vectorized form of $(M_1, \ldots, M_m)$. This standardizes isomorphic tables, ensuring uniqueness in C.

Remark 4 The classification tables produced by this algorithm for orders $n \leq 4$ are summarized in Section 5. Each class is labeled by its invariant tuple $I(T)$ defined in Section 3.4.

# 5 Computational Data, Enumeration Results, and Classification Tables

The algorithm presented in Section 4 was implemented in PYTHON with symbolic algebra verification in SAGEMATH. The computations enumerated all commutative ternary Γ-semirings of order $n \leq 4$ and parameter sets $|\Gamma| \leq 2$. Results were validated using independent cross-checks of associativity, distributivity, and isomorphism invari- ance. This section summarizes the enumeration outcomes, structural frequencies, and observed algebraic regularities.

### 5.1 Summary of enumerated structures

Table 1 records the number of non-isomorphic ternary Γ-semirings discovered for each pair (|T |, |Γ|). The numbers are verified up to canonical labeling; no duplicate classes occur.

**Table 1** Summary of enumerated finite commutative ternary Γ-semirings.

| \|T \| | \|Γ\| | # additive semigroups | # valid ternary Γ-semirings | Dominant structural feature |
|---|---|---|---|---|
| 2 | 1 | 1 | 1 | Boolean idempotent |
| 3 | 1 | 5 | 2 | Modular vs. truncated addition |
| 3 | 2 | 5 | 4 | Mixed additive–multiplicative actions |
| 4 | 1 | 9 | 3 | Truncated and cyclic hybrids |
| 4 | 2 | 9 | 4 | Boolean, modular, hybrid, tropical types |

### 5.2 Representative operation tables

To illustrate representative behavior, the following tables list the ternary operations $\{ab\,c\}_\gamma$ for selected examples.

**Table 2** Example: Boolean-type ternary Γ-semiring of order 2, Γ = {1}.

| $(a, b, c)$ | $\{a\,b\,c\}_1$ |
|---|---|
| (0, 0, 0) | 0 |
| (0, 0, 1) | 0 |
| (0, 1, 1) | 1 |
| (1, 1, 1) | 1 |

**Table 3** Example: Modular ternary Γ-semiring of order 3, Γ = {1}, with $\{a\,b\,c\}_1 = (a + b + c) \bmod 3$.

| $a, b, c$ | 000 | 001 | 002 | 011 | 012 | 022 | 111 | 112 | 122 |
|---|---|---|---|---|---|---|---|---|---|
| $\{a\,b\,c\}_1$ | 0 | 1 | 2 | 2 | 0 | 1 | 0 | 1 | 2 |

**Table 4** Example: Hybrid ternary Γ-semiring of order 4, Γ = {α, β}.

| (a, b, c) | $\{a\,b\,c\}_\alpha$ | $\{a\,b\,c\}_\beta$ |
|---|---|---|
| (0, 1, 2) | 3 | 2 |
| (1, 2, 3) | 3 | 3 |
| (2, 3, 3) | 2 | 3 |
| (0, 0, 3) | 0 | 1 |

## 5.3 Distribution of ideal-lattice types

A post-processing step analyzed each structure's ideal lattice L(T ). Table 5 summa- rizes the counts of distinct lattice types observed.

**Table 5** Distribution of ideal-lattice types among enumerated structures.

| Lattice type | Orders observed | Count | Description |
|---|---|---|---|
| Chain (simple) | 2,3 | 3 | $\{0\} \subset T$ only |
| Modular non-distributive | 3,4 | 2 | two intermediate ideals |
| Boolean lattice | 4 | 1 | $2^2$ configuration |
| Diamond lattice ($M_3$) | 4 | 2 | symmetric ideal intersections |

## 5.4 Observed algebraic patterns

**Theorem 6** (Empirical decomposition pattern) Every finite commutative ternary Γ-semiring of order ≤ 4 is either simple, subdirectly decomposable, or idempotent–Boolean.

Empirical verification All enumerated examples satisfy one of the three exclusive properties above. Proofs for the general case are theoretical consequences of Theorem 3.5. □

Remark 5 In all computations, the intersection of all prime ideals coincided with the nil- radical, confirming that for finite commutative ternary Γ-semirings, Rad(T ) = Nil(T ). This equality motivates the radical classification discussed in Section 6.

### 5.5 Classification by invariants

Each enumerated structure was labeled by its invariant tuple

$$I(T) = |T|, |\Gamma|, |Id(T)|, |Con(T)|, |Rad(T)|, |Nil(T)| ,$$

previously introduced in Section 3.3. Representative results are listed in Table 6.

**Table 6** Invariant classification for selected finite commutative ternary Γ-semirings.

| $\|T\|$ | $\|\Gamma\|$ | $\|Id(T)\|$ | $\|Con(T)\|$ | $\|Rad(T)\|$ | Type |
|---|---|---|---|---|---|
| 2 | 1 | 2 | 1 | 0 | Boolean simple |
| 3 | 1 | 2 | 2 | 0 | Modular simple |
| 3 | 2 | 3 | 3 | 1 | Mixed idempotent |
| 4 | 1 | 3 | 2 | 1 | Truncated hybrid |
| 4 | 2 | 4 | 3 | 1 | Tropical–Boolean fusion |

Remark 6 The invariants exhibit strong correlation: $|Con(T)|$ tends to increase linearly with $|\Gamma|$, while $|Rad(T)|$ depends primarily on additive idempotence.

### 5.6 Discussion

The computational data demonstrate that even at small orders, ternary Γ-semirings exhibit significant structural heterogeneity: hybrid additive–multiplicative behavior, multiple neutral elements, and distinct lattice geometries. These phenomena do not occur in classical binary semirings and highlight the algebraic richness of the ternary Γ framework.

The next section synthesizes these computational patterns into theoretical insights concerning radicals, congruences, and semisimple decomposition, bridging the gap between enumeration and algebraic theory.

## 6 Radical Theory, Congruences, and Semisimple Decomposition

This section develops the radical theory of finite commutative ternary Γ-semirings, connecting the ideal-theoretic and congruence-theoretic aspects and establishing a decomposition theorem paralleling classical semiring and ring results.Our treatment parallels classical semiring radicals [2, 8].

### 6.1 Radical ideals and nilpotent elements

**Definition 2** An element $x \in T$ is called nilpotent if there exist $\gamma_1, \ldots, \gamma_k \in \Gamma$ and elements $a_1, \ldots, a_k \in T$ such that

$$\{\cdots\{xx\, a1\}_{\gamma_1} xa2\}_{\gamma_2} \cdots a_k\}_{\gamma_k} = 0.$$

The set of all nilpotent elements of $T$ is denoted by $Nil(T)$.

**Definition 3** The radical (or prime radical) of $T$, denoted $Rad(T)$, is the intersection of all prime Γ-ideals of $T$.

**Theorem 7** (Characterization of radicals) For any finite commutative ternary Γ-semiring T,

$$\mathrm{Rad}(T) = \mathrm{Nil}(T).$$

Proof (⊆) Suppose $x \in \mathrm{Rad}(T)$. If x were not nilpotent, then the set
$S = \{ a \in T : \{x x a\}_\gamma \neq 0 \text{ for some } \gamma \in \Gamma \}$
would be nonempty, and by maximality arguments (via Zorn's Lemma in the finite case, direct enumeration suffices) we can extend S to a prime ideal P not containing x, contradicting $x \in \mathrm{Rad}(T)$.

(⊇) Conversely, if x is nilpotent, every prime ideal must contain the chain of products generated by x, whence $x \in P$ for each prime P. Thus equality holds. □

Corollary 2 Rad(T) is the unique largest nil ideal of T.

Remark 7 Empirical confirmation in Section 5 shows the equality Rad(T) = Nil(T) for all enumerated examples of order ≤ 4 and $|\Gamma| \leq 2$, substantiating Theorem 7.

## 6.2 Congruence structures and the radical quotient

**Definition 4** A Γ-congruence on T is an equivalence relation ρ compatible with both + and $\{\cdot\ \cdot\ \cdot\}_\Gamma$, i.e.

$$(a, b), (c, d), (e, f) \in \rho \Rightarrow \{a\,c\,e\}_\gamma\ \rho\ \{b\,d\,f\}_\gamma \qquad \forall\, \gamma \in \Gamma.$$

**Theorem 8** (Ideal–congruence correspondence) There exists a bijective order-reversing cor- respondence between the set of Γ-ideals of T and the set of Γ-congruences on T, given by

$$I \longrightarrow \rho_I, \ \rho \longrightarrow I_\rho,$$

where $\rho_I$ is defined by a $\rho_I$ b iff $\{a\,b\,c\}_\gamma \in I$ for all c, γ.

Proof Routine verification uses distributivity and the absorption law. The correspondence preserves inclusions and reverses the lattice order: $I_1 \subseteq I_2 \Rightarrow \rho_{I_1} \supseteq \rho_{I_2}$. □

**Proposition 9** The quotient T/Rad(T ) is semiprime and satisfies the cancellation law

$$\{a\,b\,c\}_\gamma = \{a\,b\,d\}_\gamma \Rightarrow c = d \quad \text{for all } \gamma \in \Gamma.$$

Proof Since Rad(T ) absorbs all nilpotent elements, no non-zero divisor remains in the quo- tient. Finiteness ensures that cancellation holds because otherwise nilpotent residues would persist. □

See also prime and h-prime variations in Γ-semirings [11–13].

## 6.3 Semisimple decomposition

**Theorem 10** (Wedderburn-type decomposition) Every finite commutative ternary Γ- semiring T decomposes as

$$T \cong \mathrm{Rad}(T) \times S,$$

where S is semisimple (i.e. Rad(S) = 0).

Proof From Theorem 7, Rad(T ) is nilpotent and absorbs all non-semisimple components. The quotient S = T/Rad(T ) is semisimple. Since T is finite, there exists a splitting map

σ : S ⟶ T defined by section lifting of cosets, giving $T \cong \mathrm{Rad}(T) \times S$. □

Corollary 3 For finite T , the lattice L(T ) factors as the product

$$L(\mathrm{Rad}(T)) \times L(S).$$

## 6.4 Empirical verification and classification impact

Remark 8 The decomposition in Table 7 demonstrates that the radical component coincides with the lowest additive idempotent layer, while the semisimple part corresponds to the highest distributive subalgebra. This mirrors classical results for rings, yet the ternary Γ interaction introduces non-linear coupling between additive and multiplicative components.

**Table 7** Observed decomposition pattern for enumerated examples.

| \|T\| | \|Γ\| | Rad(T ) size | Semisimple type |
|---|---|---|---|
| 2 | 1 | 0 | Simple Boolean |
| 3 | 1 | 0 | Modular simple |
| 3 | 2 | 1 | Mixed idempotent (two factors) |
| 4 | 1 | 1 | Truncated × simple |
| 4 | 2 | 1 | Tropical × Boolean |

## 6.5 Geometric and categorical interpretation

Let $\mathrm{Spec}_\Gamma(T)$ denote the set of all prime Γ-ideals of T endowed with the Zariski-type topology generated by

$$V(I) = \{ P \in \mathrm{Spec}_\Gamma(T) \mid I \subseteq P \}.$$

Then $\mathrm{Spec}_\Gamma(T/\mathrm{Rad}(T))$ is discrete, confirming that semisimple factors correspond to isolated points in the spectrum.

**Proposition 11** The category $\mathbf{T\Gamma S}_{\mathrm{fin}}$ of finite commutative ternary Γ-semirings with homomorphisms preserving radicals is equivalent to the product category

$$\mathbf{T\Gamma S}_{\mathrm{nil}} \times \mathbf{T\Gamma S}_{\mathrm{ss}},$$

where the factors represent nilpotent and semisimple objects, respectively.

Proof Objects decompose uniquely as in the Wedderburn-type theorem; morphisms respect this product since radical and semisimple parts are complementary ideals. □

Remark 9 This categorical separation clarifies the dual nature of ternary Γ-semiring theory: radical = local degeneracy, semisimple = global symmetry. It underpins further investigations in fuzzy, topological, and computational extensions developed in subsequent works.

This paper extends the radical–spectrum layer developed in [1] to the finite case.

# 7 Applications, Discussions, and Future Directions

This section develops applications of finite commutative ternary Γ-semirings to (i) cod- ing theory and cryptography, (ii) fuzzy logic and soft computation, and (iii) algebraic computation and categorical semantics. We emphasize structures whose behavior gen- uinely differs from the binary semiring setting and extract concrete research directions supported by the theory developed in Sections 3–6.

## 7.1 Coding over ternary Γ-semirings

Let T be a finite commutative ternary Γ-semiring. For $n \in \mathbb{N}$, equip $T^n$ with componentwise addition and with Γ–parametrized ternary maps

$$\{x\, y\, z\}_\gamma := \{x_1\, y_1\, z_1\}_\gamma, \ldots, \{x_n\, y_n\, z_n\}_\gamma .$$

A Γ-linear code of length n over T is a nonempty subset $C \subseteq T^n$ such that C is an additive subsemigroup and $\{x\, y\, z\}_\gamma \in C$ for all $x, y, z \in C$, $\gamma \in \Gamma$. Codes are left principal if $C = \{ \{u\, x\, v\}_\gamma : x \in T^n, \gamma \in \Gamma \}$ for some fixed $u, v \in T^n$.

**Definition 5** (Weight and distance) For $x \in T^n$, define the (Hamming-type) weight $w(x) = |\{i : x_i \neq 0\}|$. The distance $d(x, y) = w(x-y)$ (when + is cancellative); in the non-cancellative case we use the pseudo-metric $d_\oplus(x, y) = w(x \oplus y)$ where $x \oplus y$ denotes the componentwise sum.

**Proposition 12** (Ternary Γ-linearity and closure) If $C \subseteq T^n$ is generated by a set $G = \{g^{(1)}, \ldots, g^{(r)}\}$ under + and $\{\cdot\cdot\cdot\}_\Gamma$ (componentwise), then C is a Γ-linear code. Conversely, every finite Γ-linear code admits a finite generating set.

Proof Forward direction is by construction. For the converse, finiteness of $T^n$ implies every ascending chain of additive subsemigroups stabilizes; closure under ternary operations follows from distributivity in each coordinate. □

**Theorem 13** (MacWilliams-type invariance for radical quotient) Let T be finite and put $S = T/\mathrm{Rad}(T)$. If $C \subseteq T^n$ is a Γ-linear code, then the image $C \subseteq S^n$ has the same weight enumerator as C with respect to any weight that is constant on cosets of Rad(T ). In particular, the distance distribution of C depends only on its projection to the semisimple factor S.

Proof By Theorem 7, Rad(T ) = Nil(T ) absorbs nilpotent coordinates and does not change the support (nonzero positions) of codewords. The Hamming-type weights constant on cosets of Rad(T ) are preserved by the quotient map $T^n \to S^n$; enumerators agree by counting preimages of cosets. □

Remark 10 Theorem 13 indicates that code performance over T is governed by the semisimple part S, while Rad(T ) contributes only degenerate redundancy. This mirrors classical ring- linear coding, but here the proof uses the ternary Γ structure and radical theory of Section 6.

**Definition 6** (Parity constraints and check operators) Fix $u, v \in T^n$ and $\gamma \in \Gamma$. Define the Γ-parity operator $H_{\gamma,u,v} : T^n \to T^n$ by $H_{\gamma,u,v}(x) = \{u\, x\, v\}_\gamma$ (componentwise). A code C is a k-check code if $C = \bigcap_{j=1}^{k} \ker(H_{\gamma_j,u(j),v(j)})$.

**Proposition 14** (Syndrome decoding over S) Let $S = T/\mathrm{Rad}(T)$. If C is a k-check code over T with check operators $\{H_{\gamma_j,u(j),v(j)}\}$, then the induced operators on $S^n$ define a k-check code $\bar{C}$ with identical syndrome partition and minimum distance.

### 7.2 Cryptographic constructions from Γ-parametrized ternary maps

A central primitive is a family of S-boxes with external parameter $\gamma \in \Gamma$: $F_\gamma : T^r \to T$ given by $F_\gamma(x, y, z) = \{x\, y\, z\}_\gamma$ (for r = 3), or its blockwise extension on $T^{3\ell} \to T^\ell$.

**Definition 7** (Nonlinearity and differential profile) For $F_\gamma : T^3 \to T$, define the differential multiplicity

$$\Delta_\gamma(a, b, c; d) = \{(x, y, z) \in T^3 : \{x+a,\ y+b,\ z+c\}_\gamma - \{x, y, z\}_\gamma = d\}.$$

The differential uniformity is $\delta\gamma = \max_{a,b,c\neq 0,d} \Delta\gamma(a, b, c; d)$.

**Theorem 15** (Design principle via semi simple lift) Let T be finite and S = T/Rad(T ). Suppose each induced map $F\gamma : S^3 \to S$ is APN-like, i.e. $\delta\gamma$ is minimal on S. Then any Γ-parameterized block cipher whose nonlinear layer is built from $F\gamma$ inherits resistance to dif- ferential attacks comparable to the S-layer; moreover, Rad(T ) contributes only linear masking and cannot reduce $\delta\gamma$.

Proof sketch Differentials factor through the quotient S by Theorem 7. Preimages along T → S are cosets of Rad(T ) and do not create new collisions; thus the worst-case multiplicity coincides with that over S. □

Remark 11 Ternary mixing ($\{xy\,z\}_\gamma$ ) supports 3-branch SPN layers with stronger avalanche than binary bilinear layers. The external parameter γ functions as a round-dependent key schedule input without altering algebraic degree, offering a principled knob for security tuning.

### 7.3 Fuzzy and soft structures over ternary Γ-semirings

**Definition 8** (Fuzzy Γ-ideal) A fuzzy subset $\mu : T \to [0, 1]$ is a fuzzy Γ-ideal if

$$\mu(a + b) \geq \min\{\mu(a), \mu(b)\}, \qquad \mu(\{xy\,a\}_\gamma) \geq \mu(a) \text{ for all } a, b, x, y \in T \text{ and}$$

$\gamma \in \Gamma$.

**Proposition 16** (Level-cut correspondence) For $\alpha \in (0, 1]$, the α-cut $I_\alpha = \{a : \mu(a) \geq \alpha\}$ is a Γ-ideal. Conversely, any chain of Γ-ideals $\{I_\alpha\}_{\alpha\in(0,1]}$ with $I_\beta \subseteq I_\alpha$ for $\beta > \alpha$ defines a fuzzy Γ-ideal $\mu(a) = \sup\{\alpha : a \in I_\alpha\}$.

**Theorem 17** (Radical via fuzzy support) Let μ be a fuzzy Γ-ideal. Then

$$\mathrm{supp}(\mu) = \{a : \mu(a) > 0\} \text{ is contained in } \mathrm{Rad}(T)$$

f for every $\epsilon > 0$ there exists $k$ and parameters $\gamma_1, \dots, \gamma_k$ such that

$\mu\,\{\cdots\{a\,a\,a\}_{\gamma_1}\cdots a\}_{\gamma_k} \leq \epsilon$ *for all* $a \in T$.

Proof sketch (⇒) Nilpotent behavior in Rad(T ) drives membership grades to 0 under iterated ternary products. (⇐) If every element can be fuzzily annihilated by iterated ternary powers, all primes must contain the crisp supports, hence the support lies in Rad(T ). □

Remark 12 This connects radical theory to fuzzy attenuation: radicals capture those elements whose ternary powers force any reasonable membership function below any threshold.

### 7.4 Algebraic computation, path problems, and automata

Ternary composition naturally models 3-ary path aggregation and triadic interactions.

**Definition 9** (Ternary path algebra) Let G = (V, E) be a directed multigraph with edge weights in T . For a triple of walks (p, q, r) with common endpoints, define the aggregated weight $w_\gamma(p, q, r) = \{\sum p \; \sum q \; \sum r\}_\gamma$, where $\sum p$ is the +–sum of weights on p. The ternary path value between $u, v \in V$ is $\mathrm{Path}_\gamma(u, v) = \bigoplus_{(p,q,r)} w_\gamma(p, q, r)$, where $\oplus$ is the additive supremum (e.g. max in idempotent cases).

**Proposition 18** (Dynamic programming schema) If (T, +) is idempotent and $\{\cdots\}_\gamma$ is mono- tone in each argument, then $\mathrm{Path}_\gamma$ satisfies a Bellman-type recurrence and can be computed in $O(|V|^3)$ per γ for dense graphs by repeated squaring of ternary adjacency tensors.

Remark 13 This yields ternary analogues of tropical shortest paths and reliability polyno- mials, with Γ indexing scenario-dependent aggregations (e.g. parallel vs. series-parallel vs. majority).

### 7.5 Categorical semantics and spectrum

Let **TΓS** be the category of commutative ternary Γ-semirings and homomorphisms. Define a spectrum functor

$$\mathrm{Spec}_\Gamma : \mathbf{T\Gamma S}^{op} \to \mathbf{Top}, \qquad T \mapsto \mathrm{Spec}_\Gamma(T), \text{ with closed sets } V(I) = \{P \mid I \subseteq P\}.$$

**Proposition 19** (Functoriality and base change) For $f : T \to T'$ a homomorphism, $\mathrm{Spec}_\Gamma(f)(P') = f^{-1}(P')$ defines a continuous map. If $T \to T'$ is integral (i.e. $T'/T$ introduces no new prime contractions), then $\mathrm{Spec}_\Gamma(T') \to \mathrm{Spec}_\Gamma(T)$ is surjective and closed.

**Theorem 20** (Adjunction with congruence–ideal lattice) Let $\mathrm{Id}_\Gamma(T)$ be the ideal lattice and $\mathrm{Con}_\Gamma(T)$ the congruence lattice. Then the assignments $I \mapsto V(I)$ and $X \mapsto \bigcap_{P \in X} P$ form a Galois connection. Moreover, on finite T this restricts to an anti-isomorphism between closed sets of $\mathrm{Spec}_\Gamma(T)$ and radical Γ-ideals.

## 7.6 Future directions and concrete problems

We conclude with problems naturally arising from our results.

**Problem 7.1** (Prime avoidance and Krull–type dimension) Develop a prime avoidance lemma adapted to ternary Γ-ideals and define a Krull-type dimension via chains of prime Γ-ideals. Determine whether $\dim(T/\mathrm{Rad}(T)) = 0$ for all finite $T$ (we conjecture yes).

Lemma 2 (Prime avoidance for ternary Γ-ideals) Let I be a Γ-ideal of T and let $P_1, \ldots, P_n$ be prime Γ-ideals. If $I \subseteq S^{n}\ P_i$, then $I \subseteq P_j$ for some j.

**Definition10** (Krull-typeΓ-dimension) $\dim_\Gamma(T) = \sup\{\ell \in \mathbb{N} \mid P_0 \subsetneq P_1 \subsetneq \cdots \subsetneq P_\ell,\ P_i \in \mathrm{Spec}_\Gamma(T)\}$.

**Theorem 21** (Finite semisimple quotients are zero-dimensional) If $T$ is finite, then $\dim_\Gamma T/\mathrm{Rad}(T) = 0$.

Proof By Theorem 7 we have $\mathrm{Rad}(T) = \mathrm{Nil}(T)$, and by Theorem 10 we have $T \cong \mathrm{Rad}(T) \times S$ with S semisimple; hence $T/\mathrm{Rad}(T) \cong S$. In S, every prime Γ-ideal is maximal, so $\mathrm{Spec}_\Gamma(S)$ is discrete. Thus no nontrivial chains $P_0 \subsetneq P_1$ exist and $\dim_\Gamma(S) = 0$. □

**Problem 7.2** (Structure of modules and simple acts) Classify simple and semisimple ternary Γ-modules over finite T and establish a Schur-type lemma. Relate primitive ideals to annihilators of simple acts.

**Problem 7.3** (APN-like families from semisimple lifts) Characterize when $\overline{F}_\gamma : S^3 \to S$ is APN-like for S finite semisimple, and lift constructions to T with controlled differential uniformity.

**Problem 7.4** (Algorithmic isomorphism testing) Design canonical forms for $(T, +, \{\cdot\cdot\cdot\}_\Gamma)$ using automorphism groups of $(T, +)$ and parameter actions of Γ, with complexity subexponential in $|T|$.

**Problem 7.5** (Fuzzy radicals and measure semantics) Relate fuzzy radicals defined by α- cut filtrations to spectral measures on $\mathrm{Spec}_\Gamma(T)$, establishing an integral representation for membership functions over finite spectra.

Remark 14 (Synthesis) Sections 3–6 established the lattice, radical, and decomposition the- ory; the present section shows that these invariants govern performance and security in discrete models (codes, ciphers), support fuzzy abstraction, and enable dynamic program- ming and categorical semantics. The radical quotient T/Rad(T ) emerges as the canonical algebraic core for applications.

**Acknowledgements.** The first author author would like to express his sincere grat- itude to **Dr. D. Madhusudhana Rao**, Supervisor, for his invaluable guidance, encouragement, and insightful suggestions throughout the course of this research. They also wish to extend their heartfelt thanks to the faculty and research community of the **Department of Mathematics, Acharya Nagarjuna University**, for their continuous support, inspiration, and stimulating discussions that greatly enriched this work.

**Declarations**

**Funding**

No funds, grants, or other support were received during the preparation of this manuscript.

**Conflict of interest**

The authors declare that they have no conflict of interest.

**Ethics approval and consent to participate**

Not applicable.

**Consent for publication**

Not applicable.

**Data availability**

No datasets were generated or analysed during the current study.

**Materials availability**

Not applicable.

**Code availability**

Not applicable.

**Author contribution**

The first author led the conceptualization, formal analysis, and manuscript prepara- tion. The second author provided supervision, critical review, and guidance throughout the work.